\documentclass[sn-mathphys,Numbered]{sn-jnl}

\usepackage[sort&compress]{natbib}
\usepackage{graphicx}%
\usepackage{multirow}%
\usepackage{amsmath,amssymb,amsfonts,listings,mathtools,}%
\usepackage{amsthm}%
\usepackage{mathrsfs}%
\usepackage[title]{appendix}%
\usepackage{xcolor}%
\usepackage{textcomp}%
\usepackage{manyfoot}%
\usepackage{booktabs}%
\usepackage{algorithm}%
\floatname{algorithm}{Algorithm}
\usepackage{algorithmic,float}
\usepackage{listings,geometry}%
\usepackage{caption}

\theoremstyle{thmstyleone}%
\newtheorem{theorem}{Theorem}

\theoremstyle{thmstyletwo}%
\theoremstyle{thmstylethree}%
\newtheorem{definition}{Definition}%
\newtheorem{corollary}{Corollary}

\usepackage{footmisc}

\begin{document}


\title[Article Title]{Low-rank alternating direction doubling algorithm for solving large-scale continuous time algebraic Riccati equations}


\author*[1]{\fnm{Juan} \sur{Zhang}}
\email{zhangjuan@xtu.edu.cn}
\author[2]{\fnm{Wenlu} \sur{Xun}}

\affil[1]{Key Laboratory of Intelligent Computing and Information Processing
of Ministry of Education, Hunan Key Laboratory for Computation and Simulation in Science and
Engineering, School of Mathematics and Computational Science, Xiangtan University, Xiangtan,
Hunan, China}

\affil[2]{School of Mathematics and Computational Science, Xiangtan University, Xiangtan, Hunan,
China}


\abstract{This paper proposes an effective low-rank alternating direction doubling algorithm (R-ADDA) for computing numerical low-rank solutions to large-scale sparse continuous-time algebraic Riccati matrix equations. The method is based on the alternating direction doubling algorithm (ADDA), utilizing the low-rank property of matrices and employing Cholesky factorization for solving. The advantage of the new algorithm lies in computing only the $2^k$-th approximation during the iterative process, instead of every approximation. Its efficient low-rank formula saves storage space and is highly effective from a computational perspective. Finally, the effectiveness of the new algorithm is demonstrated through theoretical analysis and numerical experiments.}




\maketitle

\section{Introduction}\label{sec1}

\qquad In this paper, we consider the large-scale continuous-time algebraic Riccati equation (CARE):
\begin{equation}\label{eq1}
A^{T}X+XA+Q-XGX=0,
\end{equation}
and the complementary equation of \eqref{eq1} is
\begin{equation}\label{eq2}
AY+YA^{T}-YQY+G=0,
\end{equation}
where $A,\ G,\ Q\in\mathbb R^{n\times n},\ B\in\mathbb R^{n\times m},\ C\in\mathbb R^{p\times n}, Q=C^{T}C$, and $G=BB^{T}$, with $m, p\ll n$.

The CARE \eqref{eq1} mainly arises in the context of the quadratic optimal control problem for the following continuous-time linear time-invariant control system:
\begin{equation*}
\left\{
\begin{array}{ll}
\dot{x}(t)=Ax(t)+Bu(t),\ \ \ x(0)=x_{0},\\
y(t)=Cx(t).
\end{array}
\right.
\end{equation*}
Here, $x(t)\in\mathbb R^{n}$ is the state vector, $u(t)\in\mathbb R^{m}$ is the control vector, and $y(t)\in\mathbb R^{p}$ is the output vector. The objective of quadratic optimal control is to find a control $u(t)$ that minimizes the following function:
$$J(x_{0},u)=\frac{1}{2}\int_{0}^{+\infty}(y(t)^{T}y(t)+u(t)^{T}u(t))dt.$$
Assuming that $(A,B)$ is stabilizable and $(C,A)$ is detectable, there exists a unique optimal solution $\bar{u}$ that minimizes the functional $J(x_{0},u)$\cite{refW. M. W.}. Furthermore, this optimal solution can be determined using the feedback operator PP, such that $\bar{u}(t)=Px(t)$, where $P=B^{T}$ and $X\in\mathbb R^{n\times n}$ is the unique symmetric positive semi-definite stable solution of the CARE \eqref{eq1}.

We consider the $2n\times 2n$ Hamiltonian matrix $H$ related to the CARE \eqref{eq1}:
\begin{equation}\label{eq3}
H=\left(
        \begin{array}{cc}
          A & -G \\
          -Q & -A^{T} \\
        \end{array}
      \right),
\end{equation}
which satisfies the relation
$$HJ=-JH^{T},\ \ \ J=\left(
        \begin{array}{cc}
          0 & I_{n} \\
          -I_{n} & 0 \\
        \end{array}
      \right),$$
where $I_{n}$ denotes the identity matrix of order $n$.

Laub\cite{refA.J. Laub} proposed a numerically backward stable algorithm that applies a reordered QR algorithm\cite{refZ.J.,refJ.H.,refG.W.} to the eigenvalue problem $Hx=\lambda x$ for computing $X$. Unfortunately, the QR algorithm does not preserve the structure of the Hamiltonian matrix $H$ and the splitting of its eigenvalues. Ammar and Mehrmann\cite{refG.V.} introduced a structured-preserving algorithm that utilizes orthogonal symplectic transformations to compute a basis for the stable invariant subspace of $H$. Byers\cite{refR. Byers} presented a stable symplectic orthogonal method, but it is only applicable to systems with a single input or output. Over the past few decades, numerous iterative methods have been proposed for solving algebraic Riccati equations. The Newton method has been widely applied in the literature\cite{refL. Dieci, refS. Hammarling, refD. Kleinman,refN. Sandell}. Mehrmann and Tan\cite{refV.E.} also proposed a refinement method to correct approximate solution deficiencies. These methods require a good initial approximate solution and can be viewed as iterative refinement methods that can be combined with other direct methods. Gardiner and Laub\cite{refJ.A.J.,refR. Byers} extended the Structured Preserved Matrix Sign Function Method (MSGM)\cite{refJ. Roberts, refJ.L. Howland}. 

Additionally, a class of methods known as doubling algorithms (DA)\cite{refB.D.O.} garnered widespread interest in the 1970s and 1980s. These methods stem from fixed-point iterations derived from discrete-time algebraic Riccati equations (DARE):
$$X_{k+1}=\widehat{A}^{T}X_{k}(I+\widehat{G}X_{k})^{-1}\widehat{A}+\widehat{Q}.$$
Doubling algorithms generate sequences ${X_{2^{k}}}$ instead of sequences ${X_{k}}$. It is necessary to transform CAREs into DAREs for post-processing. However, the convergence of this algorithm has only been proven when $\widehat{A}$ is non-singular\cite{refB.D.O.}, $(\widehat{A},\ \widehat{G},\ \widehat{Q})$ is stable and detectable\cite{refM. Kimura}. The structurally preserved doubling algorithm (SDA), proposed by Guo, Lin, and Xu \cite{refX.X. Guo}, leverages the Sherman–Morrison–Woodbury formula and various iterations of sparse and low-rank representations. The resulting large-scale doubling algorithm exhibits $O(n)$ computational complexity and memory requirements per iteration and converges essentially quadratically. Furthermore, Li, Kuo, and Lin \cite{refT.X. Li} introduced a structurally preserved doubling algorithm (SDA-ls-$\varepsilon$) for solving large-scale nonsymmetric algebraic Riccati equations. They conducted a detailed error analysis of the iterative truncation effect on the approximate solution obtained by SDA.

In this paper, based on the alternating direction doubling algorithm (ADDA) proposed by Wang\cite{refW. G. Wang}, we introduce a new low-rank doubling algorithm, called R-ADDA. By transforming the Hamiltonian matrix into skew-Hermitian matrix pairs using the corresponding Cayley transformation, the low-rank alternating direction doubling algorithm inherits good convergence properties. Numerical results demonstrate that the algorithm exhibits competitiveness and effectiveness.

The structure of the remaining sections of this paper is as follows. In section \ref{sec:solve Riccati equation}, we present the iterative framework for the low-rank alternating direction doubling algorithm for solving CARE \eqref{eq1} and introduce the R-ADDA algorithm. Section \ref{sec:Convergence analysis} provides theoretical proofs of the structural preservation properties and convergence results of the new algorithm. In section \ref{sec:Numerical experiments}, we demonstrate the efficiency of the R-ADDA algorithm through numerical experiments. Finally, in section \ref{sec:Conclusions}, we present some conclusions and remarks to conclude this paper.

In this paper, we introduce some necessary symbols and terminology. Let $\mathbb{R}^{n\times m}$ denote the set of all $n\times m$ real matrices. For any matrices $A=[a_{ij}]$, $B=[b_{ij]}\in \mathbb{R}^{n\times m}$, if $a_{ij}\geq b_{ij}\ (a_{ij}>b_{ij})$ for all $i, j$, we denote $A\geq B\ (A>B)$ and define $|A|\doteq[|a_{ij}|]$. The symbols $A^{T}$ and $A^{-1}$ represent the transpose and inverse of matrix $A$, respectively. Additionally, $\rho(A)=\max_{1\leq i\leq n}\{|\lambda_{i}(A)|\}$ denotes the spectral radius of $A$, $\|A\|_{2}$ represents the 2-norm of matrix $A$, and $A\otimes I$ denotes the Kronecker product of $A$ and $I$.
\section{Low-rank alternating direction doubling algorithm}
\label{sec:solve Riccati equation}
\begin{definition}\label{def1}
For $M, L\in \mathbb{R}^{2n\times2n}$, let $M-\lambda L$ be a symplectic matrix pencil, i.e., $$MJM^{T}=LJL^{T},\ \ \ J=\left(
        \begin{array}{cc}
          0 & I_{n} \\
          -I_{n} & 0 \\
        \end{array}
      \right),$$
and define$$\mathcal{N}(M,L)=\left\{[M_{\ast},L_{\ast}]\ :\ M_{\ast},L_{\ast}\in\mathbb R^{2n\times2n},\ \text{rank}[M_{\ast},L_{\ast}]=2n,\ [M_{\ast},L_{\ast}]\left[
                                                   \begin{array}{c}
                                                     L \\
                                                     -M \\
                                                   \end{array}
                                                 \right]=0\right\}\neq\varnothing.$$
\end{definition}

\begin{definition}\label{def2}
For any given $[M_{\ast},L_{\ast}]\in\mathcal{N}(M,L)$, define
$$\widehat{M}=M_{\ast}M,\ \ \ \widehat{L}=L_{\ast}L.$$
The transformation $$M-\lambda L\rightarrow\widehat{M}-\lambda\widehat{L}$$ is called a doubling transformation.
\end{definition}

An important characteristic of this transformation is that it preserves the structure, eigenvalue spaces, and squares of eigenvalues.

Assuming $X\geq0$ is a non-negative solution to the CARE \eqref{eq1}, then the CARE \eqref{eq1} can be rewritten as
\begin{equation}\label{eq2.1}
H\left(
   \begin{array}{c}
     I \\
     X \\
   \end{array}
 \right)=\left(
   \begin{array}{c}
     I \\
     X \\
   \end{array}
 \right)R.
\end{equation}
Similarly, \eqref{eq2} yields
\begin{equation}\label{eq2.2}
H\left(
   \begin{array}{c}
     -Y \\
     I \\
   \end{array}
 \right)=\left(
   \begin{array}{c}
     -Y \\
     I \\
   \end{array}
 \right)(-S),
\end{equation}
where $$R=A-GX,\ \ S=A^{T}-QY,$$
and the matrix $H$ is the Hamiltonian matrix defined in \eqref{eq3}.

By choosing appropriate parameters $\alpha>0$ to optimize conditions for inverting certain matrices, and using Cayley transformation, we can transform \eqref{eq2.1} and \eqref{eq2.2} into the following forms
$$(H+\alpha I)\left(
   \begin{array}{c}
     I \\
     X \\
   \end{array}
 \right)(R-\alpha I)=(H-\alpha I)\left(
   \begin{array}{c}
     I \\
     X \\
   \end{array}
 \right)(R+\alpha I),$$
 $$(H+\alpha I)\left(
   \begin{array}{c}
     -Y \\
     I \\
   \end{array}
 \right)(-S-\alpha I)=(H-\alpha I)\left(
   \begin{array}{c}
     -Y \\
     I \\
   \end{array}
 \right)(-S+\alpha I).$$
 If $R-\alpha I$ and $S-\alpha I$ are nonsingular, then we have
\begin{equation}\label{eq2.3}
 \begin{array}{ll}
(H+\alpha I)\left(
   \begin{array}{c}
     I \\
     X \\
   \end{array}
 \right)=(H-\alpha I)\left(
   \begin{array}{c}
     I \\
     X \\
   \end{array}
 \right)\mathscr{C}(R;\alpha),\\
 (H+\alpha I)\left(
   \begin{array}{c}
     -Y \\
     I \\
   \end{array}
 \right)\mathscr{C}(S;\alpha)=(H-\alpha I)\left(
   \begin{array}{c}
     -Y \\
     I \\
   \end{array}
 \right),
 \end{array}
\end{equation}
where$$\mathscr{C}(R;\alpha)=(R+\alpha I)(R-\alpha I)^{-1}=(R-\alpha I)^{-1}(R+\alpha I),$$
$$\mathscr{C}(S;\alpha)=(S+\alpha I)(S-\alpha I)^{-1}=(S-\alpha I)^{-1}(S+\alpha I).$$
Assuming $A-\alpha I$ is nonsingular, and let $$A_{\alpha}=A-\alpha I,\ \ \widetilde{A}_{\alpha}=A+\alpha I,\ \ U_{\alpha}=A_{\alpha}^{T}+QA_{\alpha}^{-1}G,\ \ V_{\alpha}=A_{\alpha}+GA_{\alpha}^{-T}Q,$$
$$Z_{1}=\left(
          \begin{array}{cc}
            A_{\alpha}^{-1} & 0 \\
            QA_{\alpha}^{-1} & I \\
          \end{array}
        \right),\ \ Z_{2}=\left(
          \begin{array}{cc}
            I & 0 \\
            0 & -U_{\alpha}^{-1} \\
          \end{array}
        \right),\ \ Z_{3}=\left(
          \begin{array}{cc}
            I & A_{\alpha}^{-1}G \\
            0 & I \\
          \end{array}
        \right),$$
then we can prove that
$$M_{0}=Z_{3}Z_{2}Z_{1}(H+\alpha I)=\left(
                                           \begin{array}{cc}
                                             \widehat{A}_{0} & 0 \\
                                             -X_{0} & I \\
                                           \end{array}
                                         \right),$$
                                         $$ L_{0}=Z_{3}Z_{2}Z_{1}(H-\alpha I)=\left(
                                           \begin{array}{cc}
                                             I & Y_{0} \\
                                             0 & \widehat{A}^{T}_{0} \\
                                           \end{array}
                                         \right),$$
where \begin{equation}\label{eq2.4}
\begin{array}{lll}
\widehat{A}_{0}=I+2\alpha V_{\alpha}^{-1},\\
X_{0}=2\alpha U_{\alpha}^{-1}QA_{\alpha}^{-1},\\
Y_{0}=2\alpha A_{\alpha}^{-1}GU_{\alpha}^{-1}.
\end{array}
\end{equation}

Next, utilizing the Sherman–Morrison–Woodbury (SMW) formula
$$(A+UV^{T})^{-1}=A^{-1}-A^{-1}U(I+V^{T}A^{-1}U)^{-1}V^{T}A^{-1},$$
we can efficiently compute $U_{\alpha}^{-1}$ and $V_{\alpha}^{-1}$, then
\begin{equation}\label{eq2.5}
U_{\alpha}^{-1}=A_{\alpha}^{-T}-A_{\alpha}^{-T}C^{T}(I_{p}+CA_{\alpha}^{-1}GA_{\alpha}^{-T}C^{T})^{-1}CA_{\alpha}^{-1}GA_{\alpha}^{-T},
\end{equation}
\begin{equation}\label{eq2.6}
V_{\alpha}^{-1}=A_{\alpha}^{-1}-A_{\alpha}^{-1}B(I_{m}+B^{T}A_{\alpha}^{-T}QA_{\alpha}^{-1}B)^{-1}B^{T}A_{\alpha}^{-T}QA_{\alpha}^{-1}.
\end{equation}
We take
\begin{equation}\label{eq2.7}
\begin{array}{lll}
D_{0}=A_{\alpha}^{-T}C^{T},\ \ P_{0}=A_{\alpha}^{-1}B,\\
\Sigma_{0}=2\alpha[I_{p}-(I_{p}+D_{0}^{T}GD_{0})^{-1}D_{0}^{T}GD_{0}],\\
\Gamma_{0}=2\alpha[I_{m}-P_{0}^{T}QP_{0}(I_{m}+P_{0}^{T}QP_{0})^{-1}],
\end{array}
\end{equation}
then we can get
\begin{equation}\label{eq2.8}
X_{0}=D_{0}\Sigma_{0}D_{0}^{T},\ \ Y_{0}=P_{0}\Gamma_{0}P_{0}^{T},
\end{equation}
where $Q_{0}$ and $P_{0}$ are full-column rank matrices.

Multiplying both sides of \eqref{eq2.3} by $Z_{3}Z_{2}Z_{1}$, we get:
\begin{equation}\label{eq2.9}
\begin{array}{ll}
M_{0}\left(
   \begin{array}{c}
     I \\
     X \\
   \end{array}
 \right)=L_{0}\left(
   \begin{array}{c}
     I \\
     X \\
   \end{array}
 \right)\mathscr{C}(R;\alpha),\\
 M_{0}\left(
   \begin{array}{c}
     -Y \\
     I \\
   \end{array}
 \right)\mathscr{C}(S;\alpha)=L_{0}\left(
   \begin{array}{c}
     -Y \\
     I \\
   \end{array}
 \right).
\end{array}
\end{equation}
This iterative method can construct a pair of sequences $\{M_{k},\ L_{k}\}$, for $k=0,1,2,\cdots$, such that:
\begin{equation}\label{eq2.10}
\begin{array}{ll}
M_{k}\left(
   \begin{array}{c}
     I \\
     X \\
   \end{array}
 \right)=L_{k}\left(
   \begin{array}{c}
     I \\
     X \\
   \end{array}
 \right)[\mathscr{C}(R;\alpha)]^{2^{k}},\\
 M_{k}\left(
   \begin{array}{c}
     -Y \\
     I \\
   \end{array}
 \right)[\mathscr{C}(S;\alpha)]^{2^{k}}=L_{k}\left(
   \begin{array}{c}
     -Y \\
     I \\
   \end{array}
 \right).
\end{array}
\end{equation}
Here, $M_{k}$ and $L_{k}$ have the same form as $M_{0}$ and$ L_{0}$, such as:
$$M_{k}=\left(
         \begin{array}{cc}
          \widehat{A}_{k} & 0 \\
          -X_{k} & I \\
         \end{array}
       \right),\ \ L_{k}=\left(
                           \begin{array}{cc}
                            I & Y_{k} \\
                            0 & \widehat{A}^{T}_{k} \\
                           \end{array}
                         \right).$$

Next, we need construct $\{M_{k+1}, L_{k+1}\}$ to find suitable $\widetilde{M}, \widetilde{L}\in\mathbb R^{2n\times2n}$ such that:          $$\text{rank}((\widetilde{M},\ \widetilde{L}))=2n,\ \ (\widetilde{M},\ \widetilde{L})\left(
                                                                                  \begin{array}{c}
                                                                                    L_{k} \\
                                                                                    -M_{k} \\
                                                                                  \end{array}
                                                                                \right)=0.$$
Let $M_{k+1}=\widetilde{M}M_{k}$ and $L_{k+1}=\widetilde{L}L_{k}$, through calculations we can obtain: $\{\widetilde{M},\ \widetilde{L}\},$
$$\widetilde{M}=\left(
                  \begin{array}{cc}
                    \widehat{A}_{k}(I_{n}+Y_{k}X_{k})^{-1} & 0 \\
                    -\widehat{A}^{T}_{k}(I_{n}+X_{k}Y_{k})^{-1}X_{k} & I_{n} \\
                  \end{array}
                \right),\ \ \widetilde{L}=\left(
                                            \begin{array}{cc}
                                              I_{n} & \widehat{A}_{k}Y_{k}(I_{n}+X_{k}Y_{k})^{-1} \\
                                              0 & \widehat{A}_{k}^{T}(I_{n}+X_{k}Y_{k})^{-1} \\
                                            \end{array}
                                          \right).$$                      Then, we can get the iterative format of the Alternating Direction Doubling Algorithm (ADDA) as follows:                \begin{equation}\label{eq2.11}
\begin{array}{ll}
\widehat{A}_{k+1}=\widehat{A}_{k}(I_{n}+Y_{k}X_{k})^{-1}\widehat{A}_{k},\\
X_{k+1}=X_{k}+\widehat{A}^{T}_{k}(I_{n}+X_{k}Y_{k})^{-1}X_{k}\widehat{A}_{k},\\ Y_{k+1}=Y_{k}+\widehat{A}_{k}Y_{k}(I_{n}+X_{k}Y_{k})^{-1}\widehat{A}^{T}_{k}.
\end{array}
\end{equation}

Obviously, from the iterative format \eqref{eq2.11}, we can see that this algorithm has a computational complexity of $O(n^{3})$. By leveraging the low-rank properties of $Q=C^{T}C$ and $G=BB^{T}$, we can derive the iterative framework of the Low-Rank alternating direction doubling algorithm (R-ADDA), such that for $k=1,2,\cdots$, the R-ADDA iteration follows a recursive form:
\begin{equation}\label{eq2.12}
\begin{array}{lll}
\widehat{A}_{k}=\widehat{A}_{k-1}^{2}+\widehat{A}_{1k}\widehat{A}_{2k}^{T},\\
X_{k}=D_{k}\Sigma_{k}D_{k}^{T},\\
Y_{k}=P_{k}\Gamma_{k}P_{k}^{T},
\end{array}
\end{equation}
where $\widehat{A}_{ik}\in\mathbb R^{n\times p_{k-1}}\ (i=1,\ 2),\ D_{k}\in\mathbb R^{n\times p_{k}},\ P_{k}\in\mathbb R^{n\times m_{k}}$, and $\Sigma_{k}\in\mathbb R^{p_{k}\times p_{k}},\ \Gamma_{k}\in\mathbb R^{m_{k}\times m_{k}}$.

In the iteration, for all previous $k$, we need to store $D_{k}, \Sigma_{k}, P_{k}, \Gamma_{k}, \widehat{A}_{ik}$, and apply the Sherman-Morrison-Woodbury (SMW) formula again to obtain:
\begin{align*}
(I_{n}+Y_{k}X_{k})^{-1}&=I_{n}-Y_{k}D_{k}\Sigma_{k}(I_{p_{k}}+D_{k}^{T}Y_{k}D_{k}\Sigma_{k})^{-1}D_{k}^{T}\\
&=I_{n}-P_{k}(I_{m_{k}}+\Gamma_{k}P_{k}^{T}X_{k}P_{k})^{-1}\Gamma_{k}P_{k}^{T}X_{k},
\end{align*}
\begin{align*}
(I_{n}+X_{k}Y_{k})^{-1}&=I_{n}-D_{k}(I_{p_{k}}+\Sigma_{k}D_{k}^{T}Y_{k}D_{k})^{-1}\Sigma_{k}D_{k}^{T}Y_{k}\\
&=I_{n}-X_{k}P_{k}\Gamma_{k}(I_{m_{k}}+P_{k}^{T}X_{k}P_{k}\Gamma_{k})^{-1}P_{k}^{T},
\end{align*}
then we can obtain that
\begin{align*}
\widehat{A}_{k+1}&=\widehat{A}_{k}[I_{n}-Y_{k}D_{k}\Sigma_{k}(I_{p_{k}}+D_{k}^{T}Y_{k}D_{k}\Sigma_{k})^{-1}D_{k}^{T}]\widehat{A}_{k}\\
&=\widehat{A}_{k}[I_{n}-P_{k}(I_{m_{k}}+\Gamma_{k}P_{k}^{T}X_{k}P_{k})^{-1}\Gamma_{k}P_{k}^{T}X_{k}]\widehat{A}_{k},
\end{align*}
\begin{align*}
X_{k+1}&=X_{k}+\widehat{A}^{T}_{k}(I_{n}+X_{k}Y_{k})^{-1}X_{k}\widehat{A}_{k}\\
&=X_{k}+\widehat{A}^{T}_{k}X_{k}\widehat{A}_{k}-\widehat{A}^{T}_{k}D_{k}(I_{p_{k}}+\Sigma_{k}D_{k}^{T}Y_{k}D_{k})^{-1}\Sigma_{k}D_{k}^{T}Y_{k}X_{k}\widehat{A}_{k}\\
&=X_{k}+\widehat{A}^{T}_{k}X_{k}\widehat{A}_{k}-\widehat{A}^{T}_{k}X_{k}P_{k}\Gamma_{k}(I_{m_{k}}+P_{k}^{T}X_{k}P_{k}\Gamma_{k})^{-1}P_{k}^{T}X_{k}\widehat{A}_{k},
\end{align*}
\begin{align*}
Y_{k+1}&=Y_{k}+\widehat{A}_{k}Y_{k}(I_{n}+X_{k}Y_{k})^{-1}\widehat{A}^{T}_{k}\\
&=Y_{k}+\widehat{A}_{k}Y_{k}\widehat{A}^{T}_{k}-\widehat{A}_{k}Y_{k}D_{k}(I_{p_{k}}+\Sigma_{k}D_{k}^{T}Y_{k}D_{k})^{-1}\Sigma_{k}D_{k}^{T}Y_{k}\widehat{A}^{T}_{k}\\
&=Y_{k}+\widehat{A}_{k}Y_{k}\widehat{A}^{T}_{k}-\widehat{A}_{k}Y_{k}X_{k}P_{k}\Gamma_{k}(I_{m_{k}}+P_{k}^{T}X_{k}P_{k}\Gamma_{k})^{-1}P_{k}^{T}\widehat{A}^{T}_{k}.
\end{align*}
Here, $\oplus$ denotes the direct sum of matrices. Based on the iterative scheme above and reference \cite{refT.X. Li} , we can choose the recursive matrices as follows:
\begin{align*}
\widehat{A}_{1,k+1}&=\widehat{A}_{k}Y_{k}D_{k}\Sigma_{k}(I_{p_{k}}+D_{k}^{T}Y_{k}D_{k}\Sigma_{k})^{-1}\\
&=\widehat{A}_{k}P_{k}(I_{m_{k}}+\Gamma_{k}P_{k}^{T}X_{k}P_{k})^{-1}\Gamma_{k}P_{k}^{T}D_{k}\Sigma_{k}, \end{align*}
$$\widehat{A}_{2,k+1}=\widehat{A}_{k}^{T}D_{k},$$
$$D_{k+1}=[D_{k},\widehat{A}^{T}_{k}D_{k}],$$
$$P_{k+1}=[P_{k},\widehat{A}_{k}P_{k}],$$
\begin{align*}
\Sigma_{k+1}&=\Sigma_{k}\oplus[\Sigma_{k}-(I_{p_{k}}+\Sigma_{k}D_{k}^{T}Y_{k}D_{k})^{-1}\Sigma_{k}D_{k}^{T}Y_{k}D_{k}\Sigma_{k}]\\
&=\Sigma_{k}\oplus[\Sigma_{k}-\Sigma_{k}D_{k}^{T}P_{1k}\Gamma_{k}(I_{m_{k}}+P_{2k}^{T}X_{k}P_{k}\Gamma_{k})^{-1}P_{k}^{T}D_{k}\Sigma_{k}]\\
&=\Sigma_{k}\oplus\widetilde{\Sigma}_{k},
\end{align*}
\begin{align*}
\Gamma_{k+1}&=\Gamma_{k}\oplus[\Gamma_{k}-\Gamma_{k}P_{k}^{T}D_{k}(I_{p_{k}}+\Sigma_{k}D_{k}^{T}Y_{k}D_{k})^{-1}\Sigma_{k}D_{k}^{T}P_{k}\Gamma_{k}]\\
&=\Gamma_{k}\oplus[\Gamma_{k}-\Gamma_{k}P_{k}^{T}X_{k}P_{k}\Gamma_{k}(I_{m_{k}}+P_{k}^{T}X_{k}P_{k}\Gamma_{k})^{-1}]\\
&=\Gamma_{k}\oplus\widetilde{\Gamma}_{k}.
\end{align*}

Therefore, we can see that in the process of solving large-scale continuous-time algebraic Riccati equations using the R-ADDA iterative method, the dimensions of $D_{k}$ and $P_{k}$ grow exponentially. Moreover, the following conclusion holds:
$$\text{rank}(X_{k})\leq\text{rank}(D_{k})\leq2^{k}p,\ \ \text{rank}(Y_{k})\leq\text{rank}(P_{k})\leq2^{k}m,$$
where the number of columns of $Q_{ik}$ and $P_{ik}$ are $2^{k}p$ and $2^{k}m$, respectively. The advantage of the R-ADDA algorithm depends on the accuracy of its approximations as well as the CPU time and memory requirements.

The following is the R-ADDA algorithm for solving large-scale continuous-time algebraic Riccati equations \eqref{eq1}:
\floatname{algorithm}{Algorithm}
\renewcommand{\algorithmicrequire}{\textbf{Input:}}
\renewcommand{\algorithmicensure}{\textbf{Output:}}
\begin{algorithm}
  \caption{R-ADDA algorithm for solving CARE \eqref{eq1}}
  \label{alg1}
  \begin{algorithmic}[1]
    \REQUIRE Matrices $A\in\mathbb R^{n\times n},\ B\in\mathbb R^{n\times m},\ C\in\mathbb R^{p\times n}$, where $G=C^{T}C,\ \ Q=BB^{T}$, and parameters $\alpha$ and residual limit $\varepsilon$;
    \ENSURE $X_{k}=D_{k}\Sigma_{k}D_{k}^{T}$, such that $X_{k}\approx X$, where $X$ is the solution of CARE \eqref{eq1}.
    \STATE Compute $$A_{\alpha}=A-\alpha I,\ U_{\alpha}=A_{\alpha}^{T}+QA_{\alpha}^{-1}G,\ V_{\alpha}=A_{\alpha}+GA_{\alpha}^{-T}Q,$$ and $A_{\alpha}^{-1},\ U_{\alpha}^{-1},\ \ V_{\alpha}^{-1};$
    \STATE Set $k=0;$
    \STATE Compute $$D_{0}=A_{\alpha}^{-T}C^{T},\ \ P_{0}=A_{\alpha}^{-1}B,$$
          $$\Sigma_{0}=2\alpha[I_{p}-(I_{p}+D_{0}^{T}GD_{0})^{-1}D_{0}^{T}GD_{0}],$$
          $$\Gamma_{0}=2\alpha[I_{m}-P_{0}^{T}QP_{0}(I_{m}+P_{0}^{T}QP_{0})^{-1}];$$
    \FOR {$k=1$ until convergence}
    \STATE Compute $D_{k},\ P_{k},\ \widehat{A}_{ik}\ (i=1,2),$ and $\Sigma_{k},\ \Gamma_{k}$;
    \STATE Update $k\leftarrow k+1,$ $$X_{k+1}=D_{k+1}\Sigma_{k+1}D_{k+1}^{T};$$
    \STATE Compute $\varepsilon_{k+1}=\frac{\|A^{T}X_{k+1}-X_{k+1}A+Q-X_{k+1}GX_{k+1}\|_{2}}{\|Q\|_{2}}$ ;
            \IF {$\varepsilon_{k+1}<\varepsilon$}
                \STATE stop~(interrupt);
            \ENDIF
    \ENDFOR
  \end{algorithmic}
\end{algorithm}

Next,we consider the computational complexity of Algorithm \ref{alg1}, the operations involved include: computing $A_{\alpha}=A-\alpha I$ with a complexity of $n$; calculating $D_{0}=A_{\alpha}^{-T}C^{T},\ P_{0}=A_{\alpha}^{-1}B$ with a complexity of $4(p+m)n$; determining the complexity of $\Sigma_{0}$ and $\Gamma_{0}$ as $4(p^{2}+m^{2})n$. During the iteration, the computation complexity for $\widehat{A}_{1,k+1}$ and $\widehat{A}_{2,k+1}$ are $O((p_{k}^{3}+m_{k}^{3})$, while for $\Sigma_{k+1},\ \Gamma_{k+1}$, it is $4p_{k}m_{k}n$. Subsequently, computing $D_{k+1},\ P_{k+1}$ involves a complexity of $4(p_{k}^{2}+m_{k}^{2})n$. Therefore, the overall computational complexity sums up to $[4(p_{k}^{2}+m_{k}^{2}+p_{k}m_{k}+p^{2}+m^{2}+p+m)+1]n$.

\section{Convergence analysis }
\label{sec:Convergence analysis}
\qquad Since the R-ADDA iteration \eqref{eq2.12} is a low-rank version of the ADDA iteration \eqref{eq2.11}, we only need to present the convergence theory of the ADDA iteration \eqref{eq2.11}.

\begin{theorem}\label{thm1}\cite{refW. W. Lin}
Assuming that the matrix bundle $\widehat{M}-\lambda\widehat{L}$ is the doubling transformation of the symplectic matrix bundle $M-\lambda L$, we have the following conclusions:

(a) The matrix bundle $\widehat{M}-\lambda\widehat{L}$ is also a symplectic matrix bundle;

(b) If $M\left[
           \begin{array}{c}
             U \\
             V \\
           \end{array}
         \right]=L\left[
           \begin{array}{c}
             U \\
             V \\
           \end{array}
         \right]S$, where $U,V\in\mathbb R^{n\times m}$ and $S\in\mathbb R^{n\times m}$, then $$\widehat{M}\left[
           \begin{array}{c}
             U \\
             V \\
           \end{array}
         \right]=\widehat{L}\left[
           \begin{array}{c}
             U \\
             V \\
           \end{array}
         \right]S^{2};$$

(c) If $M-\lambda L$ has a Kronecker product canonical form $$WMZ=\left(
                                            \begin{array}{cc}
                                              J_{r} & 0 \\
                                              0 & I_{2n-r} \\
                                            \end{array}
                                          \right),\ \ \ WLZ=\left(
                                                              \begin{array}{cc}
                                                                I_{r} & 0 \\
                                                               0 & N_{2n-r} \\
                                                              \end{array}
                                                            \right),$$
where $W$ and $Z$ are non-singular, $J_{r}$ is a Jordan matrix, and $N_{2n-r}$ is a nilpotent matrix, then there exists a non-singular matrix $\widehat{W}$ such that
$$\widehat{W}\widehat{M}Z=\left(
                                            \begin{array}{cc}
                                              J_{r}^{2} & 0 \\
                                              0 & I_{2n-r} \\
                                            \end{array}
                                          \right),\ \ \ \widehat{W}\widehat{L}Z=\left(
                                                              \begin{array}{cc}
                                                                I_{r} & 0 \\
                                                                0 & N_{2n-r}^{2} \\
                                                              \end{array}
                                                            \right).$$
\end{theorem}
\begin{theorem}\label{thm2}
Assuming $X,\ Y\geq0$ are symmetric positive semi-definite solutions of equations \eqref{eq1} and \eqref{eq2} respectively, the sequences $\{\widehat{A}_{k}\},\ \{X_{k}\},\ \{Y_{k}\}$ generated by the ADDA iteration \eqref{eq2.11} satisfy:
(a) $\widehat{A}_{k}=(I_{n}+Y_{k}X)[\mathscr{C}(R;\alpha)]^{2^{k}}$;

(b) $0\leq X_{k}\leq X_{k+1}\leq X$ and $$0\leq X-X_{k}=(I_{n}+X_{k}Y)[\mathscr{C}(S;\alpha)]^{2^{k}}X[\mathscr{C}(R;\alpha)]^{2^{k}}\leq(I+XY)[\mathscr{C}(S;\alpha)]^{2^{k}}X[\mathscr{C}(R;\alpha)]^{2^{k}};$$

(c) $0\leq Y_{k}\leq Y_{k+1}\leq Y$ and $$0\leq Y-Y_{k}=(I_{n}+Y_{k}X)[\mathscr{C}(R;\alpha)]^{2^{k}}Y[\mathscr{C}(S;\alpha)]^{2^{k}}\leq(I+YX)[\mathscr{C}(R;\alpha)]^{2^{k}}Y[\mathscr{C}(S;\alpha)]^{2^{k}}.$$
\end{theorem}

\textbf{Proof.}\ \
We will prove this by induction. Firstly, we observe that $U,\ V\geq0$ implies that $I+UV$ is nonsingular, and $V(I+UV)^{-1},\ (I+UV)^{-1}U\geq0$. Given the definitions $$X_{0}=2\alpha U_{\alpha}^{-1}QA_{\alpha}^{-1},\ \ \ Y_{0}=2\alpha A_{\alpha}^{-1}GU_{\alpha}^{-1}\geq0,$$ we have $$X_{1}=X_{0}+\widehat{A}^{T}{0}(I{n}+X_{0}Y_{0})^{-1}X_{k}\widehat{A}{0}\geq X{0},$$ $$Y_{1}=Y_{0}+\widehat{A}{0}Y{0}(I_{n}+X_{0}Y_{0})^{-1}\widehat{A}^{T}{0}\geq Y{0}.$$ Since$ M_{1}-\lambda L_{1}$ is a doubling transformation of $M_{0}-\lambda L_{0}$, by equation \eqref{eq2.11}, we have:
$$\left(
    \begin{array}{cc}
      \widehat{A}_{1} & 0 \\
      -X_{1} & I \\
    \end{array}
  \right)\left(
           \begin{array}{c}
             I \\
             X \\
           \end{array}
         \right)=\left(
                   \begin{array}{cc}
                     I & Y_{1} \\
                     0 & \widehat{A}^{T}_{1} \\
                   \end{array}
                 \right)\left(
                          \begin{array}{c}
                            I \\
                            X \\
                        \end{array}
                        \right)[\mathscr{C}(R;\alpha)]^{2},$$
$$\left(
    \begin{array}{cc}
      \widehat{A}_{1} & 0 \\
      -X_{1} & I \\
    \end{array}
  \right)\left(
           \begin{array}{c}
             -Y \\
             I \\
           \end{array}
         \right)[\mathscr{C}(S;\alpha)]^{2}=\left(
                   \begin{array}{cc}
                     I & Y_{1} \\
                     0 & \widehat{A}^{T}_{1} \\
                   \end{array}
                 \right)\left(
                          \begin{array}{c}
                            -Y \\
                            I \\
                        \end{array}
                        \right),$$
then we can get
$$\widehat{A}_{1}=(I+Y_{1}X)[\mathscr{C}(R;\alpha)]^{2},\ \ \ X-X_{1}=\widehat{A}^{T}_{1}X[\mathscr{C}(R;\alpha)]^{2},$$
$$\widehat{A}^{T}_{1}=(I+X_{1}Y)[\mathscr{C}(S;\alpha)]^{2},\ \ \ Y-Y_{1}=\widehat{A}_{1}Y[\mathscr{C}(S;\alpha)]^{2},$$
which implies that $X\geq X_{1}$ and $Y\geq Y_{1}$.

Similarly, we also have
$$X-X_{1}=(I+X_{1}Y)[\mathscr{C}(S;\alpha)]^{2}X[\mathscr{C}(R;\alpha)]^{2}\leq(I+XY)[\mathscr{C}(S;\alpha)]^{2}X[\mathscr{C}(R;\alpha)]^{2},$$
$$Y-Y_{1}=(I+Y_{1}X)[\mathscr{C}(R;\alpha)]^{2}Y[\mathscr{C}(S;\alpha)]^{2}\leq(I+YX)[\mathscr{C}(R;\alpha)]^{2}Y[\mathscr{C}(S;\alpha)]^{2},$$
thus, we have proved the conclusion for $k=1$.

Next, assuming the conclusion holds for all positive integers less than or equal to $k$, we now consider the case for $k+1$. By the definition of $\widehat{A}_{k+1},\ X_{k+1},\ Y_{k+1}$ in the ADDA iteration as per equation \eqref{eq2.11}, we can similarly deduce that $$0\leq X_{k}\leq X_{k+1},\ \ \ 0\leq Y_{k}\leq Y_{k+1}.$$
Since $M_{j+1}-\lambda L_{j+1}$ is a doubling transformation of $M_{j}-\lambda L_{j},\ (j=1,1,,k)$, then by equation \eqref{eq2.11}, we have
$$\widehat{A}_{k+1}=(I+Y_{k+1}X)[\mathscr{C}(R;\alpha)]^{2^{k+1}},\ \ \ X-X_{k+1}=\widehat{A}^{T}_{k+1}X[\mathscr{C}(R;\alpha)]^{2^{k+1}},$$
$$\widehat{A}^{T}_{k+1}=(I+X_{k+1}Y)[\mathscr{C}(S;\alpha)]^{2^{k+1}},\ \ \ Y-Y_{k+1}=\widehat{A}_{k+1}Y[\mathscr{C}(S;\alpha)]^{2^{k+1}}.$$
Similarly, we can obtain
\begin{align*}
0&\leq X-X_{k+1}=(I+X_{k+1}Y)[\mathscr{C}(S;\alpha)]^{2^{k+1}}X[\mathscr{C}(R;\alpha)]^{2^{k+1}}\\
&\leq(I+XY)[\mathscr{C}(S;\alpha)]^{2^{k+1}}X[\mathscr{C}(R;\alpha)]^{2^{k+1}},
\end{align*}
\begin{align*}
0&\leq Y-Y_{k+1}=(I+Y_{k+1}X)[\mathscr{C}(R;\alpha)]^{2^{k+1}}Y[\mathscr{C}(S;\alpha)]^{2^{k+1}}\\
&\leq(I+YX)[\mathscr{C}(R;\alpha)]^{2^{k+1}}Y[\mathscr{C}(S;\alpha)]^{2^{k+1}},
\end{align*}
which proves $X\geq X_{k+1}$ and $Y\geq Y_{k+1}$. Therefore, the conclusion holds for the $k+1$case, and by the principle of mathematical induction, the theorem is proven.

Let$$W=\left[L\left[
             \begin{array}{c}
               I \\
               X \\
             \end{array}
           \right],\ M\left[
             \begin{array}{c}
               -Y \\
               I \\
             \end{array}
           \right]\right],\ \ \ Z=\left[
                                  \begin{array}{cc}
                                    I & -Y \\
                                    X & I \\
                                  \end{array}
                                \right],$$
where $M_{0}=M,\ L_{0}=L,\ X,\ Y\geq0$. From equation \eqref{eq2.9}, it follows that $W$ and $Z$ are nonsingular and satisfy
$$
W^{-1}MZ=\left(
           \begin{array}{cc}
             \mathscr{C}(R;\alpha) & 0 \\
             0 & I \\
           \end{array}
         \right),\ \ \ W^{-1}LZ=\left(
                         \begin{array}{cc}
                           I & 0 \\
                           0 & \mathscr{C}(S;\alpha) \\
                         \end{array}
                       \right).$$

Therefore, according to the spectral properties of symplectic matrix pencils, if$ \rho(\mathscr{C}(R;\alpha))<1$, then $\rho(\mathscr{C}(S;\alpha))=\rho(\mathscr{C}(R;\alpha))<1$. Additionally, if $0\leq U\leq V$, then $\|U\|_{2}\leq\|V\|_{2}$. Thus, from Theorem \ref{thm2}, we can obtain the convergence results of the ADDA algorithm as follows.
\begin{corollary}\label{cor1}
Under the assumptions of Theorem \ref{thm2}, if $\rho(\mathscr{C}(R;\alpha))<1$, then we have the following conclusions:

(a) $\|\widehat{A}_{k}\|_{2}\leq(1+\|Y\|_{2}\|X\|_{2})\|[\mathscr{C}(R;\alpha)]^{2^{k}}\|_{2}\rightarrow0\ \  (k\rightarrow\infty)$;

(b) $\|X-X_{k}\|\leq(\|X\|_{2}+\|X\|^{2}_{2}\|Y\|_{2})\|[\mathscr{C}(S;\alpha)]^{2^{k}}\|_{2}\|[\mathscr{C}(R;\alpha)]^{2^{k}}\|_{2}\rightarrow0\ \  (k\rightarrow\infty);$

(c) $ \|Y-Y_{k}\|\leq(\|Y\|_{2}+\|Y\|_{2}^{2}\|X\|_{2})\|[\mathscr{C}(R;\alpha)]^{2^{k}}\|_{2}\|[\mathscr{C}(S;\alpha)]^{2^{k}}\|_{2}\rightarrow0\ \  (k\rightarrow\infty).$
\end{corollary}

\section{Numerical experiments}
\label{sec:Numerical experiments}
\textbf{Example 1.}\ \ We take the coefficient matrix of the CARE \eqref{eq1} to be
$$A=
\left(
\begin{array}{ccccc}
    -12   &  -3   &  0   &  \cdots   &  0\\
     2 &  -12   &  -3  &   \cdots   &  0\\
     \vdots   &  \ddots  &  \ddots  &   \ddots   &  \vdots\\
     0   &  \cdots   &  2  &   -12   &  -3\\
     0   & \cdots   &  0  &   2  &   -12\\
  \end{array}
\right)_{n\times n},\ \ B=\left(
                         \begin{array}{c}
                           0.02 \\
                           0.02 \\
                           \vdots \\
                           0.02 \\
                           0.02 \\
                         \end{array}
                       \right)_{n\times1},$$
                       $$ C=\left(
                     \begin{array}{ccccc}
                       0.01, & 0.01 ,& \cdots ,& 0.01 ,& 0.01 \\
                     \end{array}
                   \right)_{1\times n}.$$

\begin{figure}[H]
\centering
    \begin{minipage}[t]{0.49\textwidth}
        \centering
        \includegraphics[width=1.1\textwidth]{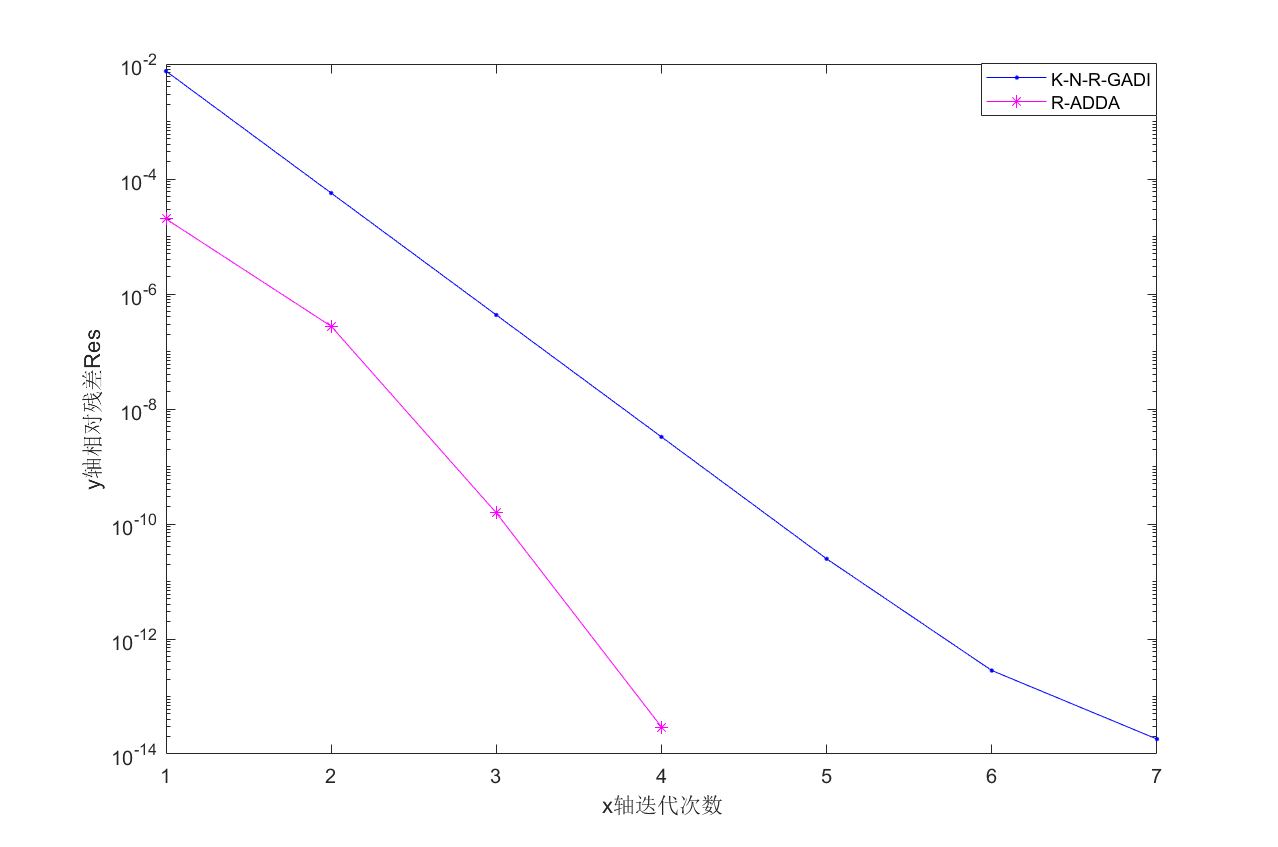}
        \caption{The residual curve of Example 1.\label{fig1}~$n=1024$}
    \end{minipage}
    \begin{minipage}[t]{0.49\textwidth}
        \centering
        \includegraphics[width=1.1\textwidth]{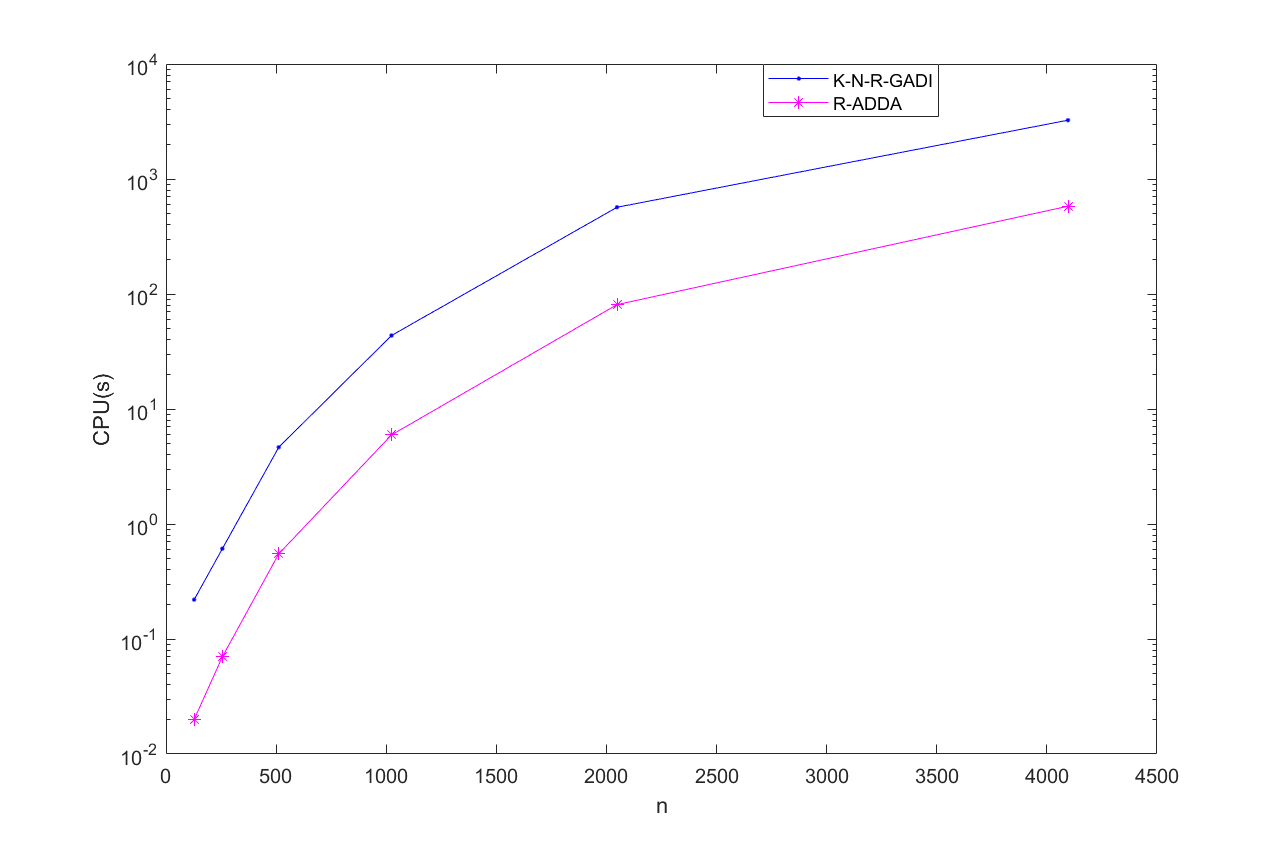}
        \caption{The time curve of Example 1.\label{fig2}}
        \end{minipage}
\end{figure}

Here, we choose the relative Residual (Res) as the iteration stopping criterion, where $$\text{Res}(X_{k})=\frac{\|A^{T}X_{k}+X_{k}A-X_{k}GX_{k}+Q\|_{2}}{\|Q\|_{2}}.$$ We use the low-rank Kleinman-Newton GADI (K-N-R-GADI) method and the R-ADDA method to solve this example, and the numerical results are shown in Table \ref{tab1}. From the table data, it is apparent that the R-ADDA iterative method is more effective in solving this example as the matrix dimension increases multiplicatively. Moreover, Figure \ref{fig1} clearly illustrates the iteration steps and the variation of relative residuals for both methods when the matrix dimension is $n=1024$, demonstrating that the convergence speed of the R-ADDA method is faster. Figure \ref{fig2} depicts the time consumption of these two iterative methods as the matrix dimension increases, highlighting the efficiency of the R-ADDA method.
\begin{table}[!htbp]
\begin{tabular}{l|l|l|l|l}
   \hline
  $n$ &algorithm & Res &IT &CPU\\
   \hline
   128& K-N-R-GADI  &  3.7511e-15 & 8 & 0.22s  \\
  \hline
   128& R-ADDA  & 6.3853e-15 & 4 & 0.02s  \\
   \hline
    256& K-N-R-GADI  &  4.3294e-14 & 7 & 0.61s  \\
   \hline
   256 & R-ADDA & 6.6167e-15 & 4 & 0.07s  \\
   \hline
    512& K-N-R-GADI  &  2.5746e-15 & 8 & 4.65s  \\
   \hline
    512 &R-ADDA & 9.1141e-15 & 4& 0.55s \\
     \hline
      1024& K-N-R-GADI  &  1.8153e-14 & 7 & 43.56s  \\
   \hline
   1024& R-ADDA  & 2.9441e-14 & 4 & 5.95s  \\
   \hline
    2048& K-N-R-GADI  &  1.6358e-14 & 7 &  567.73s  \\
   \hline
    2048& R-ADDA & 1.9252e-13 & 4 &  80.83s  \\
    \hline
    4096&K-N-R-GADI &  3.4587e-14 & 7 & 3246.6s  \\
    \hline
    4096& R-ADDA & 1.5886e-12 & 4 & 580.02s  \\
   \hline
 \end{tabular}
 \centering
 \caption{Numerical results for Example 1.\label{tab1}}
 \end{table}

\textbf{Example 2.}\ \ We take the coefficient matrix of the CARE \eqref{eq1} to be
$$A=
\left(
\begin{array}{ccccccc}
    -10   &  -3   &  -2 & 0  & 0 &  \cdots   &  0\\
     2 &  -10   &  -3  &  -2 &  0 &\cdots   &  0\\
     1 &  2   &  -10  &  -3  & -2&  \cdots   &  0\\
     \vdots   &  \ddots  &  \ddots & \ddots & \ddots &   \ddots   &  \vdots\\
     0  &  \cdots  & 1  &  2  &   -10  & -3 &  -2\\
     0 &  \cdots &0 & 1  &  2  &   -10   &  -3\\
     0   & \cdots  &0 &  0   & 1 &  2  &   -10\\
  \end{array}
\right)_{n\times n},\ \ B=\left(
                         \begin{array}{c}
                           0.005 \\
                           0.005 \\
                           \vdots \\
                           0.005 \\
                           0.005 \\
                         \end{array}
                       \right)_{n\times1},$$
                       $$C=\left(
                     \begin{array}{ccccc}
                       0.001, & 0.001 ,& \cdots ,& 0.001 ,& 0.001 \\
                     \end{array}
                   \right)_{1\times n}.$$

We use the same method as the previous example to solve this problem, and the numerical results are shown in Table \ref{tab2}. Similarly, we can conclude that compared with the K-N-R-GADI method, the R-ADDA method has a faster convergence speed and requires less iteration time, making it more effective.
\begin{table}[!htbp]
\begin{tabular}{l|l|l|l|l}
   \hline
   $n$&algorithm & Res &IT &CPU\\
   \hline
   128& K-N-R-GADI  &  4.2138e-14 & 10 & 0.34s  \\
   \hline
   128& R-ADDA  & 6.9657e-14 & 5 & 0.03s  \\
   \hline
    256& K-N-R-GADI  &  2.252e-13 & 9 & 0.85s  \\
   \hline
    256& R-ADDA & 2.5169e-13 & 5 & 0.09s  \\
   \hline
    512& K-N-R-GADI  &  1.1262e-13 & 9 & 5.31s  \\
   \hline
     512&R-ADDA & 9.5031e-13 & 5& 0.62s \\
    \hline
      1024& K-N-R-GADI  &  6.8534e-12 & 7 & 37.96s  \\
   \hline
  1024 & R-ADDA  & 3.6833e-12 & 4 & 5.16s  \\
   \hline
    2048& K-N-R-GADI  &  3.9627e-11 & 6 &  389.17s  \\
   \hline
    2048& R-ADDA & 1.4499e-11 & 4 &  61.53s  \\
   \hline
    4096&K-N-R-GADI &  1.9814e-11 & 6 & 5082.9s  \\
     \hline
   4096 & R-ADDA & 5.7516e-11 & 4 & 635.19s  \\
   \hline
 \end{tabular}
 \centering
 \caption{Numerical results for Example 2.\label{tab2}}
 \end{table}

\section{Conclusions}
\label{sec:Conclusions}
In this paper, we propose a new low-rank alternating direction doubling algorithm utilizing the low-rank properties of matrices to compute low-rank approximate solutions for large-scale algebraic Riccati equations. By combining with the ADDA algorithm, we further introduce the low-rank ADDA algorithm. Moreover, we discover that theoretically, the low-rank ADDA algorithm and the ADDA algorithm exhibit the same convergence properties. Finally, we provide numerical examples to compare the effectiveness of the low-rank Kleinman-Newton GADI algorithm and the low-rank ADDA algorithm. The results demonstrate that the low-rank ADDA algorithm is more efficient. However, like other solvers, the performance of this algorithm heavily depends on the choice of parameters, which remains a challenging issue.




\end{document}